\documentclass[12pt,a4paper,twoside]{article}

\usepackage{amsmath,amsthm,amstext,amscd,amssymb,euscript,mathrsfs}
\usepackage{calrsfs}
\usepackage{epsfig}

\newcommand{\ic}{\mathfrak{C}}

\newcommand{\Z}{\mathbb Z}
\newcommand{\R}{\mathbb R}
\newcommand{\N}{\mathbb N}

\newcommand{\E}{\mathbb E}
\newcommand{\Zd}{\mathbb Z^d}
\renewcommand{\Pr}{\mathbb P}

\renewcommand{\phi}{\varphi}

\newcommand{\La}{\ensuremath{\Lambda}}

\newcommand{\si}{\ensuremath{\sigma}}

\def\1{{\mathchoice {\rm 1\mskip-4mu l} {\rm 1\mskip-4mu l}
{\rm 1\mskip-4.5mu l} {\rm 1\mskip-5mu l}}}

\newtheorem{theorem}{{\small T}{\scriptsize HEOREM}}[section]
\newtheorem{corollary}{{\bf{\small C}{\scriptsize OROLLARY}}}[section]
\newtheorem{proposition}{{\bf{\small P}{\scriptsize ROPOSITION}}}[section]
\newtheorem{lemma}{{\bf{\small L}{\scriptsize EMMA}}}[section]
\newtheorem{remark}{{\bf{\small R}{\scriptsize EMARK}}}[section]
\newtheorem{definition}{{\bf{\small D}{\scriptsize EFINITION}}}[section]

\renewenvironment{proof}[1]
{\noindent{{\bf{\small{ P}{\scriptsize ROOF}}}.}\hspace{0.1cm} #1} {$\;\qed$\newline}

\newcommand{\beq}{\begin{eqnarray}}
\newcommand{\eeq}{\end{eqnarray}}

\newcommand{\ba}{\begin{align*}}
\newcommand{\ea}{\end{align*}}

\newcommand{\be}{\begin{equation}}
\newcommand{\ee}{\end{equation}}

\newcommand{\bl}{\begin{lemma}}
\newcommand{\el}{\end{lemma}}

\newcommand{\br}{\begin{remark}}
\newcommand{\er}{\end{remark}}

\newcommand{\bt}{\begin{theorem}}
\newcommand{\et}{\end{theorem}}

\newcommand{\bd}{\begin{definition}}
\newcommand{\ed}{\end{definition}}

\newcommand{\bp}{\begin{proposition}}
\newcommand{\ep}{\end{proposition}}

\newcommand{\bc}{\begin{corollary}}
\newcommand{\ec}{\end{corollary}}

\newcommand{\bpr}{\begin{proof}}
\newcommand{\epr}{\end{proof}}

\newcommand{\bi}{\begin{itemize}}
\newcommand{\ei}{\end{itemize}}

\newcommand{\ben}{\begin{enumerate}}
\newcommand{\een}{\end{enumerate}}

\newcommand{\caC}{{\mathscr C}}

\newcommand{\caE}{{\mathrsfs E}}

\newcommand{\caF}{{\mathcal F}}

\newcommand{\pxi}{\widehat{\mathbb{P}}_{\xi_{1}^{i-1}+_i,\xi_{1}^{i-1}-_i}}
\newcommand{\px}{\widehat{\mathbb{P}}_{X_{1}^{i-1}+_i,X_{1}^{i-1}-_i}}

\begin{document}
\title{Poincar\'e inequality for Markov random fields via disagreement percolation}  
\author{Jean-Ren\'{e} Chazottes$^{\textup{{\tiny(a)}}}$, 
Frank Redig$^{\textup{{\tiny(b)}}}$, Florian V\"{o}llering$^{\textup{{\tiny(c)}}}$\\
{\small $^{\textup{(a)}}$ Centre de Physique Th\'eorique, CNRS, \'Ecole polytechnique}\\
{\small 91128 Palaiseau Cedex, France}\\
{\small $^{\textup{(b)}}$ IMAPP, University of Nijmegen}\\{\small Heyendaalse weg
135, 6525 AJ Nijmegen, The Netherlands}\\
{\small $^{\textup{(c)}}$ Mathematisch Instituut Universiteit Leiden}\\{\small Niels Bohrweg 1, 2333 CA Leiden, The Netherlands}} 

\maketitle

\begin{abstract}
We consider Markov random fields of discrete spins on the lattice $\Zd$. 
We use a technique of coupling of conditional distributions. If under the coupling the disagreement cluster is ``sufficiently'' subcritical,
then we prove the Poincar\'e inequality. In the whole subcritical regime, we have a weak Poincar\'e inequality and corresponding polynomial upper bound for the relaxation
of the associated Glauber dynamics.

\bigskip

\end{abstract}

\section{Introduction}
Concentration inequalities is an active field of research in
probability, with applications in other areas of mathematics
such as functional analysis, geometry of metric spaces,
as well as in more applied areas such as combinatorics, optimization
and computer science \cite{ledoux}, \cite{olli}, \cite{panconesi}.

Gibbsian random fields on lattice spin systems provide 
examples of interacting random systems
having at the same time non-trivial and natural
(e.g.
Markovian) dependence
structure.
They provide 
a good class of examples where the validity of concentration
inequalities in the context of dependent random fields can be tested.

The relation between good mixing properties of Gibbs measures and
exponential relaxation to equilibrium of the associated reversible
Glauber dynamics is a thoroughly studied subject.
Well-known results in this area were obtained by 
Aizenman and Holley, \cite{AH}, Zegarli\'nski
\cite{zeg}, Stroock and Zegarli\'nski \cite{sz}, Martinelli and Olivieri \cite{MO}.
One of the main results in this area is the equivalence between the log-Sobolev
inequality (implying exponential relaxation of the dynamics
in $L^\infty$) and the Dobrushin-Shlosman complete analyticity
condition.

More recently, a direct relation between the Dobrushin uniqueness condition
and Gaussian concentration estimates was proved in \cite{kul},
and a more general relation between the existence of a coupling of a
system of conditional distributions and Gaussian and moment inequalities
in \cite{cckr}. Besides the Dobrushin uniqueness condition, disagreement
percolation technique appears here as a basic tool in constructing a good
coupling of conditional distributions. The deviation of a
function from its expectation is estimated in terms
of the sum of the squares of the maximal variation, via martingale difference
approach combined with coupling.

So far, no relation has been established between Gaussian concentration estimates
or moment estimates (such as the variance inequality) of a Gibbs measure and
relaxation properties of the associated reversible Glauber dynamics.

In this paper we show the correspondence between the existence of
a good coupling of conditional distributions and the Poincar\'{e}
inequality in the context of lattice Ising spin systems. In \cite{ccr}
this was proved in dimension one for a large class of Gibbs measures
in the uniqueness regime. The extension to higher dimension which
we deal with here (for finite-range potentials) presents new challenges.
The Poincar\'{e} inequality estimates the variance of a function
in terms of the sum of its expected quadratic variations (instead
of maximal variation). Therefore, the Poincar\'{e} inequality gives
much more information. In particular it is equivalent with
relaxation of the corresponding reversible Glauber dynamics
in $L^2$.
The Poincar\'{e} inequality is strictly weaker
than the log-Sobolev inequality. So in the complete analyticity
regime, the Poincar\'{e} inequality is satisfied. A direct
proof of the Poincar\'{e} inequality in the Dobrushin uniqueness
regime can be found in \cite{wu}.

Our result gives a direct road between ``good" coupling of
conditional distributions and the Poincar\'{e} inequality.
By good coupling we mean that if in some region of the space
we condition on two configurations
that differ only in a single point, then we can couple the unconditioned
spins such that
the set of sites
where we have a discrepancy in the coupling is small.
Small here means: behaving as a subcritical
percolation cluster, uniformly in the conditioning.
The size of this region of discrepancies can be thought of as the analogue of
the ``coupling time" for processes. In order to derive the Poincar\'{e}
inequality,
we need the existence of an exponential moment of the disagreement cluster.
which corresponds to a non-optimal
high-temperature condition (which is e.g.\ stronger than Dobrushin uniqueness, for the ferromagnetic
case). 

We want to stress however that the main message of the paper is
the direct link between coupling of conditional distributions
and the Poincar\'e inequality, rather than finding an optimal region of
$\beta$ where the inequality holds.

In case the required exponential moment of the disagreement cluster does
not exist, we still
obtain the so-called weak Poincar\'e inequality which gives
a polynomial upper bound for the relaxation of the corresponding Glauber dynamics.

Our paper is organized as follows: in section \ref{setting} we introduce the basic
ingredients and discuss coupling via disagreement percolation.
In section \ref{biboo} we prove the Poincar\'{e} inequality for small $\beta$
and $h$ close to zero, in section \ref{bambola} we treat the case $h$ large,
in section \ref{boembala} we prove the weak Poincar\'{e} inequality in
the whole subcritical regime.

\medskip

\noindent {\em Acknowledgment}. We thank Pierre Collet for fruitful discussions.

\section{Setting}\label{setting}

\subsection{Configurations}

We work in the context of Ising spin systems on a lattice, i.e., with state space
$\Omega=\{-1,+1\}^{\Zd}$ ($d\geq 2$).
Elements of $\Omega$ are denoted $\si,\eta,\xi$, and are called
spin configurations.
We fix a ``spiraling'' enumeration of $\Zd$ 
\[
 \Zd =  \{ x_1, x_2,\dots,x_n,\ldots\}.
\]
such that $x_{i+1}$ lies in the exterior boundary of $\{x_1,\ldots,x_i\}$. 
This enumeration induces an order and lattice intervals like
$$
[1, i]=\{x_k, 1\leq k\leq i\}.
$$
We use the notation $\xi_i^j$, $1\leq i\leq j\leq\infty$, for a configuration supported on the set
$\{x_k, i\leq k\leq j\}$. We denote by $\xi_{1}^{i-1}+_i$ the concatenation of $\xi_{1}^{i-1}$ with a `plus' spin
at site $x_i$. More generally, we write $\xi_V \xi_W$ for the concatenation of a configuration $\xi_V$ supported
on $V$ with a configuration $\xi_W$  supported on $W$.

\subsection{Functions}
For a function $f:\Omega\to\R$ we define the ``discrete derivative'' in the direction $\eta_x$ at the configuration $\eta$ to be
\[
\nabla_{\!\!x} f(\eta)= f(\eta^x) -f(\eta),
\]
where $\eta^x$ denotes the configuration obtained from $\eta$ by ``flipping'' the spin at site $x$, i.e., $\eta^x_y=\eta_y$ for all $y\neq x$ and $\eta^x_x=1-\eta_x$.
For a finite subset $A\subset \Zd$ we denote by $\si^A$ the configuration obtained
from $\si$ by flipping all the spins in $A$, and
$$
\nabla_{\!\!A} f(\si)= f(\si^A)-f(\si).
$$
For an enumeration $A= \{y_1,\ldots, y_n\}$ of $A$, and $x\in A$, we denote by
$A_{<x}$ the set of those elements in $A$ preceding $x$ ($x$ not included).
For the minimal element $x^*\in A$, in the chosen order of enumeration of $A$,
$A_{<x^*}=\emptyset$ by definition. 

\noindent Elementary telescoping yields the estimate
$$
|\nabla_{\!\!A} f(\si)|\leq \sum_{x\in A} \big|\nabla_{\!\!x} f(\si^{A_{<x}})\big|.
$$
Notice that if $A\subset B$ then we have the inequality
$$
\sum_{x\in A} |\nabla_{\!\!x} f(\si^{A_{<x}})|
\leq
\sum_{x\in B} |\nabla_{\!\!x} f(\si^{B_{<x}})|
$$
in an order where we enumerate $B$ by first enumerating
$A$ and then the elements of $B\backslash A$.

\noindent The variation in direction $\si_x$ is defined as
\[
\delta_x f= \sup_{\eta\in\Omega}\thinspace (f(\eta^x) -f(\eta)).
\]
The collection $\{\delta_x f:x\in\Zd\}$ is denoted by $\delta f$,
and
\[
\|\delta f\|_2^2= \sum_{x\in\Zd} \left(\delta_x f\right)^2.
\]

\subsection{Markov random fields}

Let $\mathbb{X}=\{X_x, x\in\Zd\}$ be a Markov random field of ``Ising spins'', i.e., $X_x$ takes values in $\{-1,+1\}$.
In accordance with the previous section, we use the notations $X_{1}^{i}$, $X_V$, $X_V\xi_W$, etc.

\noindent The conditional probabilities of $\mathbb{X}$ are thus given by
\be\label{ising}
\Pr\big(X_x =+1| X_{\Zd\backslash x}=\si_{\scriptscriptstyle{\Zd\backslash x}}\big)=
\frac{e^{\beta h}e^{\beta J\sum_{y\sim x}\si_y}}{2\cosh\big(\beta h+\beta J \sum_{y\sim x} \si_y\big)}\cdot
\ee
In this formula $x\sim y$ means that $x$ and $y$ are nearest neighbors, $J\in\R$ is the coupling strength and $h\geq 0$ is interpreted as a uniform magnetic field.
Without loss of generality we can assume that $|J|=1$. The case $J=1$ is the Ising ferromagnet
whereas the case $J=-1$ is the Ising anti-ferromagnet.

An easy consequence of \eqref{ising} is the following uniform bound
on the Radon-Nikodym derivative w.r.t. spin-flip:
\be\label{unibooo}
\left\|
\frac{\textup{d}\Pr^x}{\textup{d}\Pr}\right\|_\infty \leq e^{2\beta h + 4\beta d}=:e^{c}
\ee
where $\Pr^x$ denotes the image measure of $\Pr$ under
spin-flip at lattice site $x$.\newline
From the previous estimate we deduce that, for a finite subset $A\subset\Zd$,  
\be\label{unibo}
\left\|
\frac{\textup{d}\Pr^A}{\textup{d}\Pr}\right\|_\infty \leq e^{|A|(2\beta h + 4\beta d)},
\ee
where $\Pr^A$ is the image measure of $\Pr$ under simultaneous flips of all the spins in $A$.
\subsection{Glauber dynamics}
In this section we review some well-known facts about Glauber dynamics.
Much more information can be found in \cite{ligg}, chapter 3.

Given a random field $\mathbb{X}$ with distribution $\Pr$, 
the natural Glauber dynamics associated to it is
a Markovian spin-flip dynamics that flips
the spin configuration $\si$ with rate
$c(x,\si)$ at lattice site $x$.
This is the Markov process $\{\si_t:t\geq 0\}$
with
generator acting on the core of local functions given by
\be\label{genglau}
L f(\si)  = \sum_{x\in\Zd} c(x,\si) \nabla_x f(\si).
\ee
We denote by $S_t$ the associated semigroup generated by $L$, i.e.,
\[
 S_t f(\si ) = \E_\si (f(\si_t)).
\]
The rates $c(x,\si)$ are
assumed to be local, uniformly bounded away from zero and uniformly
bounded from above, i.e., there exist $0<\delta<M<\infty$ such
\be\label{boundreet}
 \delta < c(x,\si) < M.
\ee
Moreover, we assume the so-called detailed balance
relation between $c(x,\si)$ and
$\Pr$ which reads, informally,
\[
 c(x,\si) \Pr (\si) = c(x,\si^x) \Pr (\si^x).
\]
This is formally rewritten as
\be\label{libozen}
\frac{ c(x,\si)}{c(x,\si^x)}=\frac{\textup{d}\Pr^x}{\textup{d}\Pr} (\si)
\ee
i.e., the lhs of \eqref{libozen} is a (and hence the unique) continuous
(as a function of $\si$) version of the Radon-Nikodym
derivative of $\Pr$ w.r.t.\ spin-flip at site $x$ (i.e., the rhs).

Several choices for the rates are possible, one common choice
is the heat-bath dynamics where
\[
 c(x,\si) = 
\Pr\big(X_x =-\si_x| X_{\Zd\backslash x}=
\si_{\scriptscriptstyle{\Zd\backslash x}}\big).
\]
The condition \eqref{libozen} ensures that $\Pr$ is
a reversible measure for the Markov process with generator
\eqref{genglau}, i.e., the closure of $L$ is  a self-adjoint
operator on $L^2 (\Pr)$.

The Dirichlet form associated to the rates $c(x,\si)$ is
given by
\be\label{diriboy}
 \caE_c (f,f) = 2\langle f(-L)f\rangle= \sum_{x\in \Zd}\int c(x,\si) (\nabla_x f)^2 
\ \Pr (\textup{d}\si).
\ee
where $\langle \cdot \rangle$ denotes inner product in $L^2(\Pr)$.
We say that the Glauber dynamics has a spectral gap
if for all $f$ local functions with $\int f \textup{d}\Pr =0$,
\[
 \caE_c (f,f) \geq \kappa \|f\|_2^2.
\]
This implies that the $(-L)$ has simple eigenvalue zero and
that the $L^2(\Pr)$ spectrum has $\kappa$ as a lower bound.
This in turn implies the estimate
\[
 \textup{Var} (S_t f)\leq e^{-\kappa t} \|f\|_2^2 
\]
i.e., exponential relaxation to equilibrium in $L^2(\Pr)$-sense.

Defining
the quadratic form
$$
\caE (f,f) = \sum_{x\in\Zd} \int (\nabla_{\!\!x} f)^2 \textup{d}\Pr.
$$
we have by \eqref{boundreet} the estimate
\[
\delta \caE(f,f)\leq \caE_c (f,f) \leq M\caE(f,f).
\]
Hence, estimating the variance of a function in terms
of the quadratic form $\caE(f,f)$ is equivalent
with estimating the variance in terms of the Dirichlet
form \eqref{diriboy} and therefore gives relevant information
about the presence of a spectral gap and
hence $L^2$-relaxation properties of the associated Glauber
dynamics.

\subsection{Coupling of conditional probabilities}\label{coupling}

We write $\mathbb{P}_{\xi_{1}^i}$ for the conditional distribution of
$X_{[i+1,\infty)}$ given $X_{1}^i=\xi_{1}^i$.
\br
Notice that we have the same bound
\eqref{unibo} for the measure $\mathbb{P}_{\xi_{1}^i}$,
when $A\subset [1,i]^c$, uniformly in $\xi$.
\er
We denote by $\pxi$ a coupling of the distributions  
$\mathbb{P}_{\xi_{1}^{i-1}+_i}$ and $\mathbb{P}_{\xi_{1}^{i-1}-_i}$.
This coupling is a distribution of a random field
$$
\{(Y_x,Z_x), x\in [i+1,\infty)\}\quad \textup{on}\quad \big(\{-1,+1\}\times \{-1,+1\}\big)^{[i+1,\infty)}.
$$
Similarly we write $\px$. We define the random set of discrepancies 
$$
\caC_i=\{x_k : k\geq i, Y_{x_k}\neq Z_{x_k}\}.
$$
The distribution of this set depends of course on the choice of the coupling.

The coupling $\pxi$ which we will use throughout this paper is the one used in \cite{vdbm}.
For the sake of self-consistency, we explain here the construction of this
coupling.

First we pick a site $x^1_{i+1}$, with index higher than $i$, which is
a neighbor of $x_i$.
The couple $(Y_{x^1_{i+1}},Z_{x^1_{i+1}})$ is generated according
to the optimal coupling of $\mathbb{P}_{\xi_{1}^{i-1}+_i}(X_{x^1_{i+1}}=\cdot)$
and $\mathbb{P}_{\xi_{1}^{i-1}-_i}(X_{x^1_{i+1}}=\cdot)$, i.e., the coupling
that maximizes the probability of agreement.

Having generated $(Y_{x^{k-i}_{k}},Z_{x^{k-i}_{k}})$ for $i+1\leq k\leq j$, either
we choose
a new lattice point $x^{j+1-i}_{j+1}$ that has a neighbor
in the previously generated sites where $Y$ and $Z$ disagree, or if such a
point does not exist, then we choose an arbitrary neighbor higher in the
order than the previously generated sites, and generate the couple
$(Y_{x^{j+1-i}_{j+1}},Z_{x^{j+1-i}_{j+1}})$ according to the optimal coupling of 
$$\mathbb{P}_{\xi_{1}^{i-1}Y_{<}}(X_{x^{j+1-i}_{j+1}}=\cdot)
\quad
\textup{and}\quad \mathbb{P}_{\xi_{1}^{i-1}Z_{<}}(X_{x_{j+1}}=\cdot).$$
where $Y_{<}, Z_{<}$ denote the values already generated before.

By the Markov character of the random field $\mathbb{X}$, the sets of discrepancies $\caC_i$ are
almost-surely (nearest-neighbor) {\em connected}. So we can think of the $\caC_i$'s as 
{\em ``percolation clusters''} containing for sure the lattice site
$x_i$, where we have by the conditioning a disagreement.
If these clusters behave as sub-critical percolation clusters, then we say that we are in the ``good coupling regime'', see \cite{vdbm,geomaes}.
We then expect to obtain
corresponding good relaxation properties of the natural Glauber dynamics associated
to $\Pr$.
The reason to expect this is that in the entire subcritical regime for
the disagreement clusters, the corresponding Gibbs measure is unique.
In the case of the Ising model in $d=2$ it is known that in the entire
uniqueness regime we have the log-Sobolev inequality, which implies the
Poincar\'{e} inequality. It is therefore natural to expect that also in higher dimensions,
and for arbitrary Markov fields, being in the uniqueness regime implies at least
exponential
relaxation of the Glauber dynamics in $L^2$. 

\subsection{Subcritical disagreement percolation}

We suppose that, under the coupling $\pxi$, the disagreement clusters $\caC_i$ are dominated by independent subcritical site-percolation clusters, uniformly in the conditioning $\xi$. In fact, we shall need more
than subcriticality. We believe that it is an artefact of our method and that the Poincar\'e inequality
holds in the entire subcritical regime. 

We denote by $\Pr_p$ the distribution of independent site-percolation with parameter $0\leq p<1$
and by $p_c$ the corresponding critical value.
Let $\ic_i$ be the open cluster containing $x_i$. In our model \eqref{ising}, by the construction of the coupling, we have
domination by independent clusters, i.e.,
for any finite subset $A\subset\Zd$
\be\label{supsup}
\sup_i \sup_\xi \pxi \left(\caC_i\supset A \right) \leq \Pr_p(\ic_0\supset A),
\ee
with
\be\label{pising}
p=p(\beta,h)=e^{-2\beta h}\ \left(e^{4\beta d}-e^{-4\beta d}\right).
\ee
In particular,
$$
\sup_i \sup_\xi \pxi \left(|\caC_i|\geq n \right) \leq \Pr_p(|\ic|\geq n),
$$
where $\ic=\ic_0$.
Our subcriticality assumption reads as follows:
\be\label{racine}
\E_p \left(|\ic|e^{c|\ic|}\right)<\infty,
\ee 
where $c$ is defined in \eqref{unibooo}.
This condition is satisfied for $\beta$ sufficiently small or $h$ sufficiently large; see below
for the precise region of $(\beta,h)$.
\newline
By the uniform bound \eqref{supsup}, the coupling $\pxi$ can be realized in two stages.
Having generated $Y_{x_k}, Z_{x_k}$ for $k=i+1, \ldots, i+n$, we 
first generate $Y_{x_{i+n+1}}$. Then we flip an independent coin
with success probability $1-p$ (corresponding to certain agreement)
given by \eqref{pising}. Given that we have success, we
put $Z_{x_{i+n+1}}= Y_{x_{i+n+1}}$. If we do not have success, then
we possibly choose $Z_{x_{i+n+1}}= Y_{x_{i+n+1}}$ or
$Z_{x_{i+n+1}}\not= Y_{x_{i+n+1}}$ in order to obtain the correct
marginal distributions
of the coupling. The crucial point here is that the cluster
of failures (=no success), which we denote $\tilde{\caC_i}$, is a cluster that,
is {\em independent of $Y$} and contains the 
cluster of disagreement $\caC_i$. Therefore, in events that depend in a monotone
way on the cluster of disagreements $\caC_i$, we can replace
it by $\ic_i$, the cluster of failures.

\subsection{Sufficient conditions on $\beta$}

A sufficient condition for \eqref{racine} to hold is that
$$
\sum_{n=0}^\infty n\ p^n \big(2d-1\big)^n e^{cn}<\infty,
$$
where $c$ is the constant appearing in \eqref{unibooo} and $p$ 
is defined in \eqref{pising}, and where the factor 
$n\big(2d-1\big)^n$ arrizes from counting self-avoiding paths.
In turn, the above series is finite
if
$$
e^{4\beta d}-e^{-4\beta d}<\frac{e^{-4\beta d}}{2d-1},
$$
which gives
\begin{equation}\label{koko}
\beta<\frac{1}{8d} \log\Big(\frac{2d}{2d-1}\Big)\cdot
\end{equation}
Notice that this condition is independent of $h$ and 
of the sign of $J$ i.e., holds both in the ferromagnetic
and the antiferromagnetic case.
\newline
For the ferromagnetic case
$J=1$, however, the Dobrushin uniqueness condition reads
$$
2d \tanh (\beta) <1
$$
which is weaker.
See \cite{geomaes} for more details and a comparison between
uniqueness based on disagreement percolation versus
Dobrushin uniqueness.

\section{The Poincar\'e inequality and related variance inequalities }
The general idea of concentration inequalities is to give
an estimate of the probability of a deviation event $\{|f-\E(f)|> a\}$, in terms of a quantity
that measures the influence on $f$ of variations of the spin configuration
at different sites. Usually, such estimates are obtained via Chebychev's inequality, by estimating
moments of $|f-\E(f)|$, such as the variance of $f$, or higher order moments, exponential moments etc.,
in terms of a norm measuring the variability of $f$.
In this paper we concentrate on estimates of the variance.
\subsection{Uniform variance estimate}
The semi-norm
\[
 \|\delta f\|^2_2 = \sum_{x\in \Zd} (\delta_x f)^2
\]
measures the influence of spin-flips on $f$ in a uniform way, i.e., for each
$x$ the worst influence is computed.

The first inequality measures the variance in terms
of $\|\delta f\|_2^2$.
\bd\label{devroyeineq}
We say that a random field $\mathbb{X}$ satisfies the uniform
variance inequality 
if there exists
$C>0$, such for all $f:\Omega\to\R$, $f\in L^2 (\Pr)$,  we have
\be\label{devroyebamboe}
\E( (f-\E(f))^2)\leq C \|\delta f\|_2^2
\ee
\ed
The uniform variance inequality estimates
the variance in terms of the rather ``rough'' norm $\|\delta f\|_2^2$.
Surprisingly,
it is still a powerful inequality with many useful applications, such as
almost-sure central limit theorems, convergence of the
empirical distribution in a strong (Kantorovich) distance,
etc. See \cite{colletchaz} for a list of applications.

Examples where the uniform variance inequality is satisfied include
high-temperature Gibbsian random fields (where
it follows from the much stronger log-Sobolev inequality)
and plus phase 
of the Ising model
at low enough temperatures, see \cite{cckr}.

\subsection{Poincar\'{e} inequality}

The quadratic form
$$
\caE (f,f) = \sum_{x\in\Zd} \int (\nabla_{\!\!x} f)^2 \textup{d}\Pr
$$
measures the influence of spin-flips on $f$, taking into account the
distribution of the spin-configuration, i.e., 
large differences between $f(\si^x)$ and $f(\si)$
are weighted less if they correspond to exceptional configurations (in the sense
of the measure $\Pr$).
We have the obvious inequality $\caE(f,f)\leq \|\delta f\|_2^2$, therefore, estimating
the variance in terms of $\caE(f,f)$ is clearly better, and, as we will see in
examples below, this difference can be substantial.

\bd\label{poincaredef}
We say that the random field $\mathbb{X}$ satisfies the Poincar\'e inequality if there exists a 
constant $C_{\!\scriptscriptstyle{P}}>0$ such that for all $f\in L^2(\Pr)$ 
\begin{equation}\label{poincare}
\int \big(f-\E(f)\big)^2 \textup{d}\Pr\leq C_{\!\scriptscriptstyle{P}}\ \caE (f,f)\ .
\end{equation}
\ed

The Poincar\'{e} inequality is strictly stronger than the uniform
variance inequality. Moreover, contrary to the uniform variance estimate,
the Poincar\'{e} inequality gives exponentially fast decay to equilibrium
for the associated Glauber dynamics
in $L^2 (\Pr)$. Indeed, 
\eqref{poincare} implies
\[
\textup{ Var}(f) \leq \frac{1}{\delta}C_{\!\scriptscriptstyle{P}} \caE_c(f,f)
= 2\langle f,(-L)f \rangle
\]
from which one easily sees that $(-L)$ has a spectral
gap in $L^2(\Pr)$ of at least $\kappa=2\delta/C_{\!\scriptscriptstyle{P}}$, which implies
the relaxation estimate
\[
\textup{ Var} (S_t f) \leq e^{-\kappa t} \|f\|_2^2
\]
\subsection{Weak Poincar\'{e} inequality}

Finally, the variance can be estimated in terms of a combination
of $\caE(f,f)$ and another term $\Phi (f)$, where
$\Phi$ is homogeneous of degree 2, i.e.,
$\Phi(\lambda f)=\lambda^2\Phi(f)$. Examples are
$\Phi(f)=\|f\|_\infty^2$, or $\Phi (f)=\|\delta f\|_2^2$.
The idea here is that
if the Poincar\'{e} inequality does not hold, it can be due to
``bad
events'' which have relatively small probability (e.g.\ large disagreement
clusters). The idea is then to estimate the variance
by $\caE(f,f)$ on the good configurations and
by $\Phi(f)$ on the bad configurations.
This leads to the weak  Poincar\'{e} inequality,
initially introduced by R\"{o}ckner and Wang
\cite{rw}. This inequality contains enough information to conclude 
relaxation properties of the associated Glauber dynamics,
but now with $\textup{Var} (S_t f)$ estimated with a stronger norm than the $L^2(\Pr)$-norm.
\bd\label{weakpoincarebamboe}
The measure $\Pr$ satisfies the weak Poincar\'e inequality 
if there exists a decreasing function
$\alpha: (0,\infty)\to (0,\infty)$ such that
for all bounded measurable functions
$f:\Omega\to\R$
we have, for all $r>0$
$$
\int \big(f-\E(f)\big)^2 \textup{d}\Pr\leq \alpha(r)\ \caE (f,f) + r\Phi(f)\ .
$$
\ed
If we have
\be\label{kontra}
\Phi( S_t f)\leq \Phi(f)
\ee
i.e., if $S_t$ contracts $\Phi(\cdot)$, then we obtain
a relaxation estimate for the dynamics from the weak Poincar\'{e} inequality. More precisely,
in that case, for bounded measurable
functions $f$  with $\int f \textup{d}\Pr=0$, the weak Poincar\'e inequality implies the estimate
$$
\textup{Var} (S_t f) \leq \xi(t)\left(\|f\|_2^2 + \Phi(f)\right)
$$
where $\xi(t)\to 0 $ as $t\to\infty$ is determined by $\alpha$:
\[
\xi (t)= \inf\left\{r>0: -\frac{1}{\delta}\alpha(r)\log r\leq 2t\right\},\;t>0.
\]
where $\delta>0$ is the lower-bound on the spin-flip rates.
In the case when $\alpha(r)\leq C r^{-\kappa}$ for $C,\kappa>0$,
we get $\xi(t)\leq \big(1+\frac{1}{\kappa}\big)^{1+\frac{1}{\kappa}} \left(\frac{2t\delta}{C}\right)^{-\frac{1}{\kappa}}$.
We refer the reader to \cite{rw} for more background and details.

\subsection{Examples}
Here we illustrate with some simple examples that the Poincar\'{e} inequality
is much stronger than the uniform variance inequality. The examples
are representants of a whole class of functions for which the effect
of spin-flip is only ``typically small'', which gives a good estimate
of $\caE(f,f)$, but where the uniform variation $\delta_i f$ is
always of order one.

Let $d=1$ and $\Pr$ be a translation invariant probability
measure
on configurations $\si\in \Omega = \{ -1,+1\}^{\Z}$
such that there exists $0<\theta<1$ with
\be\label{soliboze}
 \Pr (\si_1=\alpha_1,\ldots,\si_n=\alpha_n )\leq \theta^n
\ee
for all $n\in\N$, $\alpha_1,\ldots,\alpha_n\in \{-1,1\}$.
Examples of such $\Pr$ are translation-invariant Gibbs measures.

Consider for $n\in \N, k<n$
\[
 f_k (\si_1,\ldots, \si_n) = \left|\left\{i\in \{ 1,\ldots,n-k\}:
\si_i= \si_{i+1}=\cdots=\si_{i+k}=+1\right\}\right|
\]
i.e., the number of lattice intervals of size $k$, contained in $[1,n]$
and filled with plus spins. 

We have
\[
 \nabla_r f_k (\si) = \sum_{j\in [1,n-k]: r\in[j,j+k]} 
\left(\1\{\si_r=-1\}-\1\{\si_r=+1\}\right)\prod_{i\in[j,j+k], i\not=r} \1\{\si_i=+1\} 
\]
which gives
\[
\int (\nabla_r f_k)^2\textup{d}\Pr \leq 2k\theta^k
\]
and hence
\[
 \caE(f_k,f_k) \leq 2k(n-k) \theta^k.
\]
Therefore, if $\Pr$ satisfies the Poincar\'{e} inequality
(e.g.\ for a large
class of Gibbs measures in one dimension in the uniqueness regime,
\cite{ccr}) then
\[
 \textup{Var} (f_k)\leq C_{\!\scriptscriptstyle{P}} 2k (n-k) \theta^k
\]
Choosing now
$k = c\log (n)$, and putting $\theta= e^{-\alpha}$ we find
that
\[
\textup{Var}( f_{c\log n})\leq 2 c\log (n) (n-c\log (n)) n^{-\alpha c}.
\]

Hence if $\alpha c>1$, 
$\textup{Var}( f_{c\log n})$ goes to zero as $n\to\infty$.
It is immediate from \eqref{soliboze} that $\alpha>c$ the first moment 
$\E( f_{c\log (n)})$ converges
to zero as $n\to\infty$. Therefore, $\alpha c>1$ implies that
$f_{c\log n}$ converges to zero in $L^2(\Pr)$ (and hence in probability)
as $n\to\infty$.

On the other hand, it is clear that $\delta_i (f)= 1$ for all
$i=1,\ldots,n$, therefore the uniform variance estimate
gives $\textup{Var}(f_k) \leq Cn$, which is not useful here.

One can consider similar quantities like the number
of clusters of size $k$ of plus-spins, the number
of self-overlaps of size $k$, etc. Such quantities
will have small $\caE(f,f)$ (for measures satisfying \eqref{soliboze})
and large $\|\delta f\|_2^2$.

\section{Poincar\'{e} inequality for the case $h=0$ }\label{biboo}

We start with the following result.

\bt\label{colibri}
Consider the Markov random defined in \eqref{ising} with $h=0$. For $\beta$ chosen such that
$$
\E_p \left(|\ic|e^{c|\ic|}\right)<\infty,
$$
the Poincar\'e inequality \eqref{poincare} holds.
\et

In section \ref{bambola} below (Theorem \ref{colobro}), we will give a complementary result
which covers the case of large $\beta$ and (correspondingly) large $h$.

\bpr
The proof is divided in four steps.\\
\noindent {\bf Step 1} (Martingale decomposition).\newline
Let $f:\Omega\to\R$ be a bounded measurable function. Define
$$
\Delta_i=\Delta_i(X_1^i)= \E(f| \caF_i)-\E(f| \caF_{i-1})
$$
where $\caF_i$ is the sigma-field generated by $\{X_{x_k} : 1\leq k\leq i\}$ for $i\geq 1$
and where $\caF_0$ is the trivial sigma-field $\{\emptyset,\Omega\}$. Then we have
$$
\text{Var}(f) = \sum_{i\in\N} \E(\Delta_i^2).
$$
\noindent {\bf Step 2} (Coupling representation of $\Delta_i$)\newline
We have (using that spins can take only two values)
\begin{align}
\nonumber
|\Delta_i | &= \left|\int \textup{d}\Pr_{X_1^{i-1}}(\xi_i)
\int \textup{d}\widehat{\mathbb{P}}_{X_1^i,X_1^{i-1}\xi_i} (\si_{i+1}^\infty,\eta_{i+1}^\infty)
\left( f(X_1^{i-1} X_i\si_{i+1}^\infty) - f(X_1^{i-1} \xi_i\eta_{i+1}^\infty) \right)
\right|\\
\nonumber
& \leq \int \left| f(X_1^{i-1} +_i\si_{i+1}^\infty) - f(X_1^{i-1} -_i\eta_{i+1}^\infty) \right|
\textup{d}\px (\si_{i+1}^\infty,\eta_{i+1}^\infty)\\
\nonumber
&= \int \left| f(X_1^{i-1} +_i\si_{i+1}^\infty) - f(X_1^{i-1} -_i\eta_{i+1}^\infty) \right|
\textup{d}\px  (\si_{i}^\infty,\eta_{i}^\infty)\\
\nonumber
&= \sum_{A\ni x_i}\int \textup{d}\px  (\si_{i}^\infty,\eta_{i}^\infty)\; \times \\
\label{floflo}
&\qquad \1\{\caC_i=A\}\left| f(X_1^{i-1}\eta_{A}\si_{(A\cup[1,i-1])^c}) - f(X_1^{i-1} \si_{A}\si_{(A\cup[1,i-1])^c}) \right|,
\end{align}
where $\widehat{\mathbb{P}}_{X_1^{i-1}+_i,X_1^{i-1}-_i}$ is the coupling of conditional
probabilities defined in subsection \ref{coupling}. Notice that the sum over $A$ runs over 
{\em finite} connected subsets of $\Zd$
containing $x_i$ since $\caC_i$ is dominated by a subcritical percolation cluster. 

In the sequel, we simply write $\si_V\xi_W\eta$ for $\si_V\xi\eta_{(V\cup W)^c}$ to alleviate
notations.
\bigskip

\noindent {\bf Step 3} (Telescoping and domination by independent clusters).

\noindent Start again from \eqref{floflo} and telescope the disagreement cluster:
\beq
|\Delta_i|
&\leq &
\int \left|\nabla_{\caC_i} f(X_1^{i-1}\sigma^{\caC_i})\right| \ \textup{d}\px (\si, \eta)
\nonumber\\
&\leq &
\int \sum_{x\in\caC_i} \big|\nabla_x f(X_1^{i-1}\sigma^{(\caC_i)_{<x}})\big| \ \textup{d}\px (\si, \eta)
\nonumber\\
&\leq &
\int \sum_{x\in\tilde{\caC}_i} \big|\nabla_x f(X_1^{i-1}\sigma^{(\tilde{\caC}_i)_{<x}})\big| \ \textup{d}\px (\si, \eta)
\nonumber\\
&=&
\widetilde{\E}\int\sum_{x\in\tilde{\caC_i}}
\big|\nabla_x f(X_1^{i-1}\sigma^{(\tilde{\caC_i})_{<x}})\big|
\ \textup{d}\Pr_{X^{i-1}_1+_i} (\si)
\nonumber\\
&=& \sum_{A\ni x_i} \sum_{x\in A} \Pr_p (\ic_i =A)\int \big|\nabla_x f(X_1^{i-1}\sigma^{A_{<x}})\big|
 \ \textup{d}\Pr_{X^{i-1}_1+_i} (\si).
 \nonumber
\eeq
In the third inequality the expectation is over the ``failure cluster"
$\tilde{\caC_i}$ only, which is independent of $\si$. This independence gives 
the factorization in the last equality, by decomposing over the realization
of this cluster (which is finite with $\px$ probability one under the subcriticality assumption).
\bigskip

\noindent {\bf Step 4} (Change of measure).

\noindent Using now the bound \eqref{unibo} and the remark in the beginning
of subsection \ref{coupling}, we further estimate, using
$$
|\Delta_i|
\leq
\sum_{A\ni x_i} \sum_{x\in A} \Pr_p (\ic_i =A)\ e^{c|A|}
\int |\nabla_{\!\!x} f(X_1^{i-1}+_i\si)| \
\textup{d}\Pr_{X^{i-1}_1+_i} (\si)
$$
where $c$ is defined in \eqref{unibooo}.\newline
Define the finite number (by the subcriticality assumption \eqref{racine})
$$
K:= \sum_{A\ni 0} |A| \ \Pr_p (\ic=A) \ e^{c|A|} = \E_p \big(|\ic| e^{c|\ic|}\big).
$$
Then, using the elementary inequality
\be\label{element}
\left(\sum_{k} a_k b_k\right)^2\leq \sum_{k} a_k \sum_k a_k b_k^2
\ee
for $a_k, b_k \geq 0$, we obtain
\beq
\sum_{i\in\N} \E (\Delta_i^2)
&\leq &
Ke^{2c} \sum_{i\in\N}\sum_{A\ni x_i}\sum_{x\in A} e^{c|A|}  \Pr_p (\ic_i = A)
\int (\nabla_{\!\!x} f)^2 \textup{d}\Pr
\nonumber\\
&=&
K^2e^{2c} \caE (f,f),
\nonumber
\eeq
where the extra factor $e^c$ arises from removing the plus in
the conditioning in $\Pr_{X_1^{i-1} +_i}$.
This finishes the proof of Theorem \ref{colibri}
\epr
\section{ Non-zero magnetic field}\label{bambola}

In this section we show how to prove the Poincar\'{e} inequality under a subcricality condition different from Theorem \ref{colibri}.
It is strictly worse in the case $h=0$ (since it uses Cauchy-Schwarz to seperate the realization of the disagreement
cluster from the gradient of $f$) but can be used for $\beta$ large and $h$ large, where the condition \eqref{racine}
fails.

\bt\label{colobro}
Suppose that $p$ given in \eqref{pising} is such that
\be\label{racine2}
\sum_n n(2d-1)^n e^{c'n} \Pr_p (|\ic|\geq n)^{1/2} <\infty,
\ee
where 
\be\label{dising}
c'=4\beta d.
\ee
Then the Poincar\'{e} inequality holds. 
\et

For \eqref{racine2} to hold, it is sufficient that
$$
(2d-1) p^{\frac12} e^{c'} <1
$$
which gives
$$
(2d-1)^2 e^{-2\beta h}(e^{12\beta d}-e^{4\beta d})<1.
$$
This is satisfied for $\beta$ small enough or $h$ large enough.

\bpr
The telescoping and coupling steps are the same
as in the proof of Theorem 1. So we arrive at
$$
|\Delta_i|  \leq 
\sum_{A\ni x_i}\sum_{x\in A} \int \textup{d}\px 
(\si_{i}^\infty,\eta_{i}^\infty)  \1\{\caC_i=A\} \big|\nabla_{\!\!x} f(X_1^{i-1}\si_{A_{<x}}\eta)\big|.
$$
Now we use Cauchy-Schwarz inequality to obtain
\begin{align}
\nonumber
|\Delta_i|  & \leq 
\sum_{A\ni x_i}\sum_{x\in A}  \left(\px\big(\caC_i=A\big)\right)^{1/2} \ \times \\
& \qquad \qquad\left(
\int \textup{d}\px (\si_{i}^\infty,\eta_{i}^\infty)  \Big(\nabla_{\!\!x} f(X_1^{i-1}\si_{A_{<x}}\eta)\Big)^2
\right)^{1/2}.
\label{CS}
\end{align}

\bigskip

\noindent {\bf Step 4} (Change of measure). In the r.h.s. of \eqref{CS} we integrate over
the ``composite'' configuration $\si_{A_{<x}}\eta$ under the coupling
$\px$. To recover the measure $\Pr$ (see
later) we need to replace $\si_{A_{<x}}$ by $\eta_{A_{<x}}$. The cost of this replacement
is independent of $h$ and
is estimated in the following lemma where $\pxi$ is the coupling introduced above. 

\bl
Let $A$ be a finite subset of $\Zd$ containing $x_i$
and let $x\in A$. Let $\Pr_1$ be the distribution of $Z_{A_{<x}}Y_{(A_{<x})^c}$ and
$\Pr_2$ be the distribution of $\{Y_x,x\in\Zd\}$. Then $\Pr_1$ is absolutely continuous
with respect to $\Pr_2$ and 
$$
\left\|\frac{\textup{d}\Pr_1}{\textup{d}\Pr_2}\right\|_\infty \leq e^{c'|A|}
$$
where $c'$ is defined in \eqref{dising}.
\el
\bpr
Let $\La\subset \Zd$ finite, large enough to contain $A$. We have
by construction of the coupling $\pxi$  (see subsection \ref{coupling}):
\begin{align*}
& \frac{\pxi\big(Z_{A_{<x}}=\si_{\!\scriptscriptstyle{A_{<x}}},Y_{\La\backslash A_{<x}}=
\eta_{\scriptscriptstyle{\La\backslash A_{<x}}}\big)}
{\pxi\big(Y_{A_{<x}}=\si_{\!\scriptscriptstyle{A_{<x}}},Y_{\La\backslash A_{<x}}=
\eta_{\scriptscriptstyle{\La\backslash A_{<x}}}\big)}\\
&=
\sum_{\zeta_{A_{<x}}}
\Pr_{\xi_1^{i-1}-_i}(\si_{\!\scriptscriptstyle{A_{<x}}})\
\Pr_{\xi_1^{i-1}+_i\zeta_{A_{<x}}}(\eta_{\scriptscriptstyle{\La\backslash A_{<x}}})\times \\
& \quad\quad\quad\quad\frac{\pxi\big(Z_{A_{<x}}=\zeta_{\scriptscriptstyle{A_{<x}}}\big|Y_{A_{<x}}=
\si_{\!\scriptscriptstyle{A_{<x}}}\big)}
{\Pr_{\xi_1^{i-1}+_i\si_{\!\scriptscriptstyle{A_{<x}}}}\big(\eta_{\scriptscriptstyle{\La\backslash A_{<x}}}\big)}\\
&\leq \sup_{\zeta} \ \frac{\Pr_{\xi_1^{i-1}+_i\zeta_{\scriptscriptstyle{A_{<x}}}}
(\eta_{\scriptscriptstyle{\La\backslash A_{<x}}})}
{\Pr_{\xi_1^{i-1}+_i\si_{\!\scriptscriptstyle{A_{<x}}}}\big(\eta_{\scriptscriptstyle{\La\backslash A_{<x}}}\big)}\\
& \leq e^{c' |\partial A_{<x}|}\leq e^{c' |A|}.
\end{align*}
We conclude by letting $\La\uparrow\Zd$.
\epr

Returning to \eqref{CS} and using the preceding lemma we get
\begin{align}
\nonumber
|\Delta_i|  & \leq 
\sum_{A\ni x_i}\sum_{x\in A}  \left(\px\big(\caC_i=A\big)\right)^{1/2} \ e^{c'|A|} \ \times \\
\nonumber
& \qquad \qquad  \qquad \qquad \left(
\int \textup{d}\Pr_{X_1^{i-1}-_i}(\eta)  \big(\nabla_{\!\!x} f(X_1^{i-1}\eta)\big)^2
\right)^{1/2}\\
\nonumber
& \leq e^c \ \sum_{A\ni x_i}\sum_{x\in A}  \left(\px\big(\caC_i=A\big)\right)^{1/2} \ e^{c'|A|} \ \times \\
\label{fluflu}
& \qquad \qquad  \qquad \qquad \left(
\int \textup{d}\Pr_{X_1^{i}}(\eta)  \big(\nabla_{\!\!x} f(X_1^{i}\eta)\big)^2
\right)^{1/2},
\end{align}
where for the second inequality we used that, under the measure $\Pr$, the cost
of  flip at a single site is bounded by $e^c$ (see \eqref{unibooo}).

\bigskip

\noindent {\bf Step 5} (Domination by independent clusters).
Using \eqref{supsup} we get from \eqref{fluflu}
\begin{align}
\nonumber
|\Delta_i|  
& \leq e^c \ \sum_{A\ni x_i}\sum_{x\in A}  \big(\Pr_p(|\ic|\geq |A|)\big)^{1/2} \ e^{c'|A|} \ \times \\
\label{flifli}
& \qquad \qquad  \qquad \qquad \left(
\int \textup{d}\Pr_{X_1^{i}}(\eta)  \big(\nabla_{\!\!x} f(X_1^{i}\eta)\big)^2
\right)^{1/2}.
\end{align}
Now let 
$$
K'=\sum_{A\ni x_i}\sum_{x\in A} \Pr_p(|\ic|\geq |A|)^{1/2}\ e^{c'|A|} =
\sum_{A\ni 0} \ |A| \ \Pr_p(|\ic|\geq |A|)^{1/2}\ e^{c'|A|}. 
$$
By assumption \eqref{racine2} $K'$ is finite. Using once more
the elementary inequality \eqref{element}
we deduce from \eqref{flifli} that
\begin{align*}
\sum_i \E(\Delta_i^2)
& \leq e^{2c} K'\ \sum_i \sum_{A\ni x_i}\sum_{x\in A}  \Pr_{\!p}\big(|\ic|\geq |A|\big)^{1/2} \ e^{c'|A|} \ 
\int\big(\nabla_{\!\!x} f\big)^2 \textup{d}\Pr\\
& =
e^{2c} K'\ \sum_x \left(\int \big(\nabla_{\!\!x} f\big)^2\textup{d}\Pr\right)\ \sum_{A\ni x}\ |A| \
\Pr_p\big(|\ic|\geq |A|\big)^{1/2} \ e^{c'|A|}\\
&= C_{\!\scriptscriptstyle{P}}\ \sum_x \int \big(\nabla_{\!\!x} f\big)^2\ \textup{d}\Pr
\end{align*}
where
$$
C_{\!\scriptscriptstyle{P}}:= e^{2c} K'^2.
$$
This finishes the proof of Theorem \ref{racine2}.
\epr

\section{Weak Poincar\'{e} inequality}\label{boembala}

If the assumption \eqref{racine} fails, but $p<p_c$
(where $p_c$ denotes the critical value for independent
site percolation) then we are still in the uniqueness
regime (i.e., the conditional probabilities \eqref{ising} admit
a unique Gibbs measure) and expect suitable decay properties
of the Glauber dynamics. 

We show that in this regime the weak Poincar\'e inequality holds,
which gives polynomial relaxation to equilibrium.

\bt
Suppose that $p$ \textup{(}defined in \eqref{pising}\textup{)}
satisfies $p<p_c$. Then the weak Poincar\'{e}
inequality is satisfied.
Moreover, there exists $C,\kappa>0$ such that
\[
\alpha (r)\leq  C r^{-\kappa}.
\]
As a consequence, 
\[
\textup{Var}(S_t f)\leq  \left(1+\frac{1}{\kappa}\right)^{1+\frac{1}{\kappa}} \left(\frac{2t\delta}{C}\right)^{-\frac{1}{\kappa}} \left(\|f\|_2^2 +4\|f\|_\infty^2\right) 
\]
where
$\delta$ is defined in \eqref{boundreet}.
\et
\bpr
The proof follows the lines of the proof of Theorem 1,
so we sketch where we start to deviate from it:
In the estimation of the variance, the contribution involving $\|f\|_\infty^2$ will arise
by cutting the cluster of disagreement at some order
of magnitude $N$. 

The sum in \eqref{racine} is now possibly infinite, so we define
$$
K_N=\sum_{n=0}^N  n\ e^{cn} \ \Pr_p\left(|\ic|\geq n \right).
$$
Following the line of proof of Theorem \ref{colibri},
we follow the change of measure road for realizations
of the cluster $\ic_i =A$ of cardinality less than or equal to $N$, and for
$A$ with $|A|>N$ we use the uniform estimate
\[
\sup_{\eta} |f(\eta^A)- f(\eta)|\leq \sum_{x\in A} \delta_x f\leq 2|A|\|f\|_\infty
\]
Next estimate, using Jensen and the elementary inequality \eqref{element},
\beq
&&\sum_{i\in\N} \Big(\sum_{A\ni x_i, |A|>N} \Pr_p (\ic_i=A)
\sum_{x\in A}(\delta_x f)\Big)^2
\nonumber\\
&\leq &
4\left(\E_p ( |\ic|^2 \1\{|\ic|>N\})\right)^2\|f\|_\infty^2.
\nonumber
\eeq
This gives the inequality
$$
\textup{Var}(f) \leq 2e^c K_N^2 \caE (f,f) + 8 \left(\E_p ( |\ic|^2 \1\{|\ic|>N\})\right)^2\|f\|_\infty^2.
$$
The constant in front of $\caE (f,f)$ blows up at most exponentially in $N$,  i.e., we have
the estimate
\[
2e^c K_N^2 \leq C_1 e^{a N}
\] 
where $C_1,a$ are strictly positive and $(\beta,h)$-dependent.
The constant in front of  $\|f\|_\infty^2$ is exponentially
small in the whole subcritical regime, by the exponential
decay of the cluster size, \cite{grim} i.e., we have the estimate
\[
2\left(\E_p ( |\ic|^2 \1\{|\ic|>N\})\right)^2\leq C_2 e^{-bN}
\]
where $C_2,b$ are strictly positive  and $(\beta,h)$-dependent. Therefore we can take
$$
\alpha(r)\leq C_1 \left(\frac{r}{C_2}\right)^{-\frac{a}{b}}
$$
and $\kappa=a/b$.  
\epr
  
\end{document}